\documentclass[12pt]{article}
\usepackage{amsbsy,amsfonts}

\newtheorem{theo}{Theorem}[section]
\newtheorem{lemma}[theo]{Lemma}

\newtheorem{coll}[theo]{Corollary}
\newtheorem{remark}[theo]{Remark}

\newcommand{\be}{\begin{equation}}
\newcommand{\ee}{\end{equation}}
\newcommand{\ba}{\begin{array}}
\newcommand{\ea}{\end{array}}

\newcommand{\m}{\frak m}
\newcommand{\dsum}{\displaystyle \sum}
\newcommand{\dprod}{\displaystyle \prod}

\begin{document}

\title{Links of symbolic powers of prime ideals}
\author{ Hsin-Ju Wang \thanks{e-mail:  hjwang@math.ccu.edu.tw} \\
Department of Mathematics, National Chung Cheng University\\
          Chiayi 621, Taiwan}

\date{}

\maketitle

\begin{abstract}
~\\ In this paper, we prove the following. Let $(R, \m )$ be a
Cohen-Macaulay local ring of dimension $d\geq 2$. Suppose that
either $R$ is not regular or $R$ is regular with $d\geq 3$. Let
$t\geq 2$ be a positive integer. If $\{\alpha_1, \dots,
\alpha_d\}$ is a regular sequence contained in $\m^t$, then
$$(\alpha_1, \dots, \alpha_d): \m^t\subseteq  \m^t.$$ This result
gives an affirmative answer to a conjecture raised by Polini and
Ulrich.
\end{abstract}

\section{Introduction}

Let $R$ be a Noetherian ring. Two proper ideals $I$ and $L$ of
height $g$ are said to be {\it directly linked}, $I\sim L$, if
there exists a regular sequence ${\bf z}=z_1, \dots, z_g\subseteq
I\cap L$ such that $I=({\bf z}):L$ and $L=({\bf z}):I$. If $R$ is
Cohen-Macaulay and $I$ is an unmixed ideal of height $g$ with
$R_P$ Gorenstein for every minimal prime $P$ of $I$, then one can
always produce a link $L\sim I$ by choosing a regular sequence
${\bf z}=z_1, \dots, z_g\subset I$ and setting $L=({\bf z}):I$
(\cite{ps}). The linkage class of a proper ideal $I$ is the set of
all ideals $K$ obtained by a sequence of links $I=L_0\sim L_1\sim
\cdots \sim L_n=L$, where $n\geq 0$. The algebraic foundations of
linkage theory were established in \cite{hu} and \cite{ps}.
\par In \cite{pu}, Polini and Ulrich investigated when an ideal is
the unique maximal element of its linkage class, in the sense that
it contains every ideal of the class. They showed that if $(R,
\m)$ is a Gorenstein local ring of dimension $d\geq 2$, with
$d\geq 3$ if $R$ is regular, then every ideal in the linkage class
of $\m^t$ is contained in $\m^t$ provided that $gr_{\m}(R)$ is
Cohen-Macaulay, or $R$ is a complete intersection, or
$ecodim~R\leq 3$, or $t\leq 3$. They pointed out that $\m^t$ will
be the maximal element of its linkage class if the following
conjecture holds. \vskip 0.05in ~~\\  {\bf Conjecture~1:} Let $(R,
\m)$ be a local Gorenstein ring of dimension $d\geq 2$, with
$d\geq 3$ if $R$ is regular, and let $t\geq 1$ be an integer; then
every
ideal directly linked to $\m^t$ is contained in $\m^t$. \\
In the same paper, they also showed that Conjecture~1 is
equivalent to the following question: Is a two-dimensional local
Gorenstein ring $(R, \m)$ regular if and only if there exists an
integer $t\geq 1$ and a regular sequence ${\bf \alpha}=\alpha_1,
\alpha_2$ contained in $\m^t$ so that $({\bf \alpha}): \m^t$ is
not contained in $\m^t$? \vskip 0.1in  In this paper, we settle
the conjecture in the affirmative by showing the following:

\begin{theo}\label{theo4}
Let $(R, \m )$ be a Cohen-Macaulay local ring of dimension $d\geq
2$, with $d\geq 3$ if $R$ is regular. Let $t\geq 2$ be a positive
integer. If $\{\alpha_1, \dots, \alpha_d\}$ is a regular sequence
contained in $\m^t$, then
$$(\alpha_1, \dots, \alpha_d): \m^t\subseteq  \m^t.$$
\end{theo}
\par Let $R$ be a Noetherian ring. Recall that an ideal $I$ of
$R$ of height $g$ is said to be equi-multiple of reduction number
1 if there exists a $g$ generated ideal $J$ contained in $I$ so
that $I^2=JI$. Ideals that are equi-multiple of reduction number 1
enjoy some nice properties: For example, it is shown in \cite{va}
that if $R$ is a Cohen-Macaulay local ring and $I$ is an
equi-multiple Cohen-Macaulay ideal of reduction number 1, then the
associated graded ring $gr_I(R)$ of $I$ is Cohen-Macaulay.
Therefore, it is of great interest to search for such classes of
ideals. In \cite{cp1} and \cite{cpv}, Corso, Polini and
Vasconcelos proved that the links of prime ideals in
Cohen-Macaulay rings are equi-multiple ideals of reduction number
1 under certain assumptions. Moreover, Corso and Polini \cite{cp2}
addressed the following conjecture. \vskip 0.05in ~~\\
{\bf Conjecture~2:} Let $(R, \m)$ be a Cohen-Macaulay local ring
and $P$ a prime ideal of height $g\geq 2$. Let $J=(z_1, \dots,
z_g)$ be an ideal generated by a regular sequence contained in
$P^{(s)}$, where $s$ is a positive integer $\geq 2$. For any $k=1,
\dots, s$ set $I_k=J:P^{(k)}$. Then
$$I^2_k=JI_k$$ if one of the following two conditions holds:
\begin{description} \item{($IL_1$)} $R_P$ is not a regular local
ring; \item{($IL_2$)} $R_P$ is a regular local ring and two of the
$z_i$'s lie in $P^{(s+1)}$.
\end{description}

With the help of Theorem~\ref{theo4} and Theorem~\ref{theo8}, we
are able to give a positive answer to Conjecture~2.

\begin{theo} \label{theo5}
Let $R$ be a Noetherian ring. Let $P$ be a prime ideal of height
$g\geq 2$, $J=(z_1, \dots, z_g)$ an ideal generated by a regular
sequence contained in $P^{(k)}$, where $k\geq 2$ is a positive
integer. Set $I_k=J:P^{(k)}$. Then $I_k^2=I_kJ$ if one of the
following conditions holds:
\begin{description} \item{(i)} $R_P$ is not a regular local ring.
\item{(ii)} $R_P$ is a regular local ring and $g\geq 3$.
\item{(iii)} $R_P$ is a regular local ring, $g=2$, and $z_i\in
P^{(k+1)}$ for every $i$.
\end{description}
\end{theo}
\section{Main Theory}
In this section, we study Conjecture~1 stated in section 1. We are
able to prove it in a more general form.
\begin{theo} \label{theo6}
Let $(R, \m )$ be a Cohen-Macaulay local ring of dimension $d\geq
2$. Let $I$ be an ideal of $R$ of height 2 and $\tilde{I}=\sum
((\alpha, \beta): I)$, where the sum is taken over all regular
sequences $\{\alpha, \beta\}$ in $I$. Let $t$ be a positive
integer. If $\{\alpha, \beta\}$ is a regular sequence contained in
$I^t$, then
$$(\alpha, \beta): I^t\subseteq I^{t-1}\tilde{I}.$$
\end{theo}

\begin{remark}
\emph{(i) $I\subseteq \tilde{I}\subseteq R$.\\ (ii) Conjecture~1
follows if we set $d=2$ and $I=\m$ in Theorem~\ref{theo6} and
observe that $\tilde{I}=I=\m$ if $R$ is not regular.\\ (iii) One
can see from the proof of Theorem~\ref{theo6} that the assumption
that $R$ is Cohen-Macaulay can be removed if the grade of $I$ is
2. }
\end{remark}
To show Theorem~\ref{theo6}, we begin with the following lemmas.
\begin{lemma} \label{lem1}
Let $(R, \m)$ be a Noetherian local ring. Let $\{x_1, \dots,
x_{t-1}, y_1, \dots, y_{t-1}\}\subseteq \m$ such that $\{x_i,
y_j\}$ is a regular sequence $\forall i, j$, where $t\geq 2$. Then
$$\ba{rl} & \cap_{i=1}^{t-1} (x_1\cdots x_i, y_1\cdots y_{t-i})\\
= & <\{x_1\cdots x_{i-1}y_1\cdots y_{t-i}~|~1\leq i\leq t\}>
\ea.$$
\end{lemma}
\begin{proof}
It is enough to show that $$\ba{rl} & \cap_{i=1}^k (x_1\cdots x_i,
y_1\cdots y_{t-i})\\ = & (x_1\cdots x_k)+<\{x_1\cdots
x_{i-1}y_1\cdots y_{t-i}~|~1\leq i\leq k\}> \ea$$ for $1\leq k\leq
t-1$. By induction on $k$, suppose that the equality holds for
$k$. Then $$\ba{rl} & \cap_{i=1}^{k+1} (x_1\cdots x_i, y_1\cdots
y_{t-i})\\ = & ((x_1\cdots x_k)+<\{x_1\cdots x_{i-1}y_1\cdots
y_{t-i}~|~1\leq i\leq k\}>)\cap (x_1\cdots x_{k+1}, y_1\cdots
y_{t-k-1})\\  = & (x_1\cdots x_k)\cap (x_1\cdots x_{k+1},
y_1\cdots y_{t-k-1})+ <\{x_1\cdots x_{i-1}y_1\cdots
y_{t-i}~|~1\leq i\leq k\}> \\ = &  (x_1\cdots
x_{k+1})+<\{x_1\cdots x_{i-1}y_1\cdots y_{t-i}~|~1\leq i\leq
k+1\}>. \ea$$

\end{proof}

\begin{lemma} \label{lem2}
Let $(R, \m )$ be a Cohen-Macaulay local ring of dimension $d\geq
2$. Let $I$ be an ideal of $R$ of height 2 and $\tilde{I}=\sum
((\alpha, \beta): I)$, where the sum is taken over all regular
sequences $\{\alpha, \beta\}$ in $I$. Let $t\geq 1$ and $\{x_1,
\dots, x_t, y_1, \dots, y_t\}\subseteq I$ such that $\{x_i, y_j\}$
is a regular sequence $\forall i, j$. Then
$$(x_1\cdots x_t, y_1\cdots y_t): I^t\subseteq I^{t-1}\tilde{I}.$$
\end{lemma}

\begin{proof}
We may assume that $t\geq 2$. First observe that $$(x_1\cdots x_t,
y_1\cdots y_t): I^t\subseteq \cap_{i=1}^{t-1} (x_1\cdots x_i,
y_1\cdots y_{t-i}).$$ This is because that
$$\ba{rl} (x_1\cdots x_t, y_1\cdots y_t):
I^t & \subseteq (x_1\cdots x_t, y_1\cdots y_t): x_{i+1}\cdots
x_ty_{t-i+1}\cdots y_t \\ & =(x_1\cdots x_i, y_1\cdots y_{t-i})
\ea$$ for $1\leq i\leq t-1$. By Lemma~\ref{lem1}, every element in
$(\prod_i x_i, \prod_i y_i): I^t$ is of the form
$$\sum_{i=1}^ta_ix_1\cdots x_{i-1}y_1\cdots y_{t-i}.$$ Therefore,
to show the lemma, it is enough to show that $a_i\in \tilde{I}$
for every $i$. For this, let $z\in I$; then
$$(\sum_{j=1}^ta_jx_1\cdots x_{j-1}y_1\cdots
y_{t-j})zx_{i+1}\cdots x_ty_{t-i+2}\cdots y_t\in (x_1\cdots x_t,
y_1\cdots y_t).$$ Observe that $$(a_jx_1\cdots x_{j-1}y_1\cdots
y_{t-j})zx_{i+1}\cdots x_ty_{t-i+2}\cdots y_t\in (x_1\cdots x_t,
y_1\cdots y_t)$$ if $j\neq i$. Therefore
$$a_izx_1\cdots \hat{x}_i \cdots x_ty_1\cdots \hat{y}_{t-i+1}\cdots y_t\in
(x_1\cdots x_t, y_1\cdots y_t).$$ It follows that $a_iz\in (x_i,
y_{t-i+1})$, which implies that $a_i\in (x_i, y_{t-i+1}):
I\subseteq \tilde{I}$.
\end{proof}

\vskip 0.2in ~~\\{\bf Proof of Theorem~\ref{theo6}:} We first
prove the following statement: If $y_i\in I~\forall i$ and
$\{\alpha, y_1\cdots y_t\}$ is a regular sequence contained in
$I^t$, then
$$(\alpha, y_1\cdots y_t): I^t\subseteq I^{t-1}\tilde{I}.$$ Observe that the following
set generates $I^t$: $$S_t=\{z_1\cdots z_t~|~ z_j\in I,~\{y_i,
z_j\}~is~a~regular~sequence~\forall i, j\}.$$ Let $a\in (\alpha,
y_1\cdots y_t): I^t$. Let $z_1\cdots z_t\in S_t$; then there are
elements $u(z), v(z)\in R$ such that \be \label{eq1} az_1\cdots
z_t=u(z)\alpha+v(z)(y_1\cdots y_t).\ee Therefore if $z'_1\cdots
z'_t$ is any element in $S_t$ then \be  \label{eq2} az'_1\cdots
z'_t=u(z')\alpha+v(z')(y_1\cdots y_t).\ee From (\ref{eq1}),
(\ref{eq2}) and the fact that $\{\alpha, y_1\cdots y_t\}$ is a
regular sequence, we see that $u(z)z'_1\cdots z'_t\in (z_1\cdots
z_t, y_1\cdots y_t)$. Since $z'_1\cdots z'_t$ is an arbitrary
element in $S_t$,
$$ u(z)\in (z_1\cdots z_t, y_1\cdots y_t): I^t\subseteq I^{t-1}\tilde{I}$$
by Lemma~\ref{lem2}. Since $z_1\cdots z_t$ is an arbitrary element
in $S_t$ and $\alpha \in I^t$, there are elements $u, v\in R$ with
$u\in I^{t-1}\tilde{I}$ such that $a\alpha=u\alpha+v(y_1\cdots
y_t)$, which implies that $a-u\in (y_1\cdots y_t)$. It follows
that $a\in I^{t-1}\tilde{I}$ and this completes the proof of the statement.\\
In general, if $\{\alpha, \beta\}$ is a regular sequence contained
in $I^t$, then we can use the same argument the above to conclude
that every element in $(\alpha, \beta): I^t$ is in
$I^{t-1}\tilde{I}$.

\begin{coll}
Let $(R, \m )$ be a Cohen-Macaulay local ring of dimension $d\geq
2$. Suppose that either $R$ is not regular or $R$ is regular with
$d\geq 3$. Let $t\geq 2$ be a positive integer. If $\{\alpha_1,
\dots, \alpha_d\}$ is a regular sequence contained in $\m^t$, then
$$(\alpha_1, \dots, \alpha_d): \m^t\subseteq  \m^t.$$
\label{coll1}
\end{coll}

\begin{proof}
We prove by induction on $d$. The case $d=2$ follows from
Theorem~\ref{theo6}.
\\ Assume that $d\geq 3$. Let $a\in (\alpha_1, \dots, \alpha_d):
\m^t$. Since $R/(\alpha_d)$ is Cohen-Macaulay and is not regular,
$a\in \m^t~(mod~(\alpha_d))$ by induction. It follows that $a\in
\m^t$ as $\alpha_d\in \m^t$.
\end{proof}

~~\\ The following result will be used in section 3.
\begin{coll}
Let $(R, \m )$ be a Cohen-Macaulay local ring of dimension $d\geq
2$. Suppose that either $R$ is not regular or $R$ is regular with
$d\geq 3$. Let $t\geq 2$ be a positive integer and $\{\alpha_1,
\dots, \alpha_d\}$ be a regular sequence contained in $\m^t$. Then
for every $a\in (\alpha_1, \dots, \alpha_d): \m^t$ and every
element $f\in \m^t$, if $af=\sum_{i=1}^d u_i\alpha_i$ then $u_i\in
\m^t$ for every $i$. \label{coll2}
\end{coll}

~~\\ Corollary~\ref{coll2} follows from the next lemma.
\begin{lemma}\label{lem3}
Let $(R, \m )$ be a Cohen-Macaulay local ring of dimension $d\geq
2$. Let $t\geq 2$ be a positive integer. Assume that for any
regular sequence $\{\alpha_1, \dots, \alpha_d\}$ contained in
$\m^t$, $(\alpha_1, \dots, \alpha_d): \m^t\subseteq \m^t$. Then
for every $a\in (\alpha_1, \dots, \alpha_d): \m^t$ and every
element $f\in \m^t$, if $af=\sum_{i=1}^d u_i\alpha_i$ then $u_i\in
\m^t$ for every $i$.
\end{lemma}

\begin{proof}
It is easy to see that the following set generates the ideal
$\m^t$:
$$S_t=\{z_1\cdots z_t~|~ z_j\in \m,~\{\alpha_1, \dots, \hat{\alpha}_i, \dots,
\alpha_d, z_j\}~is~a~regular~sequence~\forall i, j\}.$$ Let
$z_1\cdots z_t\in S_t$; then there are $u_i(z)\in R$ such that
$$az_1\cdots z_t=\sum_{i=1}^d u_i(z)\alpha_i.$$ From the proof of
Theorem~\ref{theo6}, we see that $u_i(z)\in (\alpha_1, \dots,
\hat{\alpha}_i, \dots, \alpha_d, z_1\cdots z_t): \m^t$ for every
$i$. Therefore $u_i(z)\in \m^t$ by assumption. This shows that if
$f\in \m^t$ then $af$ can be expressed as $\sum_{i=1}^d
u_i\alpha_i$ for some $u_i\in \m^t$ for every $i$. If
$af=\sum_{i=1}^d u'_i\alpha_i$ for some $u'_i\in R$, then
$u_i-u'_i\in (\alpha_1, \dots, \hat{\alpha}_i, \dots,
\alpha_d)\subseteq \m^t$, it follows that $u'_i\in \m^t$ for every
$i$.
\end{proof}

\section{Ideals of reduction number 1}

The goal of this section is to study the following conjecture: \vskip 0.03in ~~\\
{\bf Conjecture:} Let $(R, \m)$ be a Cohen-Macaulay local ring and
$P$ a prime ideal of height $g\geq 2$. Let $J=(z_1, \dots, z_g)$
be an ideal generated by a regular sequence contained in
$P^{(s)}$, where $s$ is a positive integer $\geq 2$. For any $k=1,
\dots, s$ set $I_k=J:P^{(k)}$. Then
$$I^2_k=JI_k$$ if one of the following two conditions holds:
\begin{description} \item{($IL_1$)} $R_P$ is not a regular local
ring; \item{($IL_2$)} $R_P$ is a regular local ring and two of the
$z_i$'s lie in $P^{(s+1)}$.
\end{description}

There are some partial answers to this conjecture in \cite{cp2}
and \cite{pu}. The following result solves the conjecture
completely.

\begin{theo} \label{theo7}
Let $R$ be a Noetherian ring. Let $P$ be a prime ideal of height
$g\geq 2$, $J=(z_1, \dots, z_g)$ an ideal generated by a regular
sequence contained in $P^{(k)}$, where $k\geq 2$ is a positive
integer. Set $I_k=J:P^{(k)}$. Then $I_k^2=I_kJ$ if one of the
following conditions holds:
\begin{description} \item{(i)} $R_P$ is not a regular local ring.
\item{(ii)} $R_P$ is a regular local ring and $g\geq 3$.
\item{(iii)} $R_P$ is a regular local ring, $g=2$, and $z_i\in
P^{(k+1)}$ for every $i$.
\end{description}
\end{theo}
To show Theorem~\ref{theo7}, we need the following result.

\begin{theo}\label{theo8}
Let $(R, \m)$ be a Cohen-Macaulay local ring of dimension $d\geq
2$. Let $k\geq 2$ be a positive integer. If ${\bf z}=z_1, \dots,
z_d$ is a regular sequence contained in $\m^{k+1}$, then $$ ({\bf
z}): \m^k\subseteq (\m^{2k+1}: \m^k)$$ if depth$~gr_{\m}(R)\geq
2$.
\end{theo}
\begin{proof}
It is not difficult to see that the following set generates the
ideal $\m^k$.
$$ S_k=\{\dprod_{i=1}^k y_i~|~y_i\in \m,  \{y^{\ast}_i, y^{\ast}_j\}~is ~
a~regular~sequence~in~gr_{\m}(R)\}.$$ Let $a\in ({\bf z}): \m^k$
and $f=\dprod_{i=1}^k y_i \in S_k$; then \be \label{eq3}
af=\dsum_{i=1}^d u_iz_i\ee for some $u_i\in R$. To show the
theorem it is enough to show that $u_i\in \m^k$ as $z_i\in
\m^{k+1}$ for every $i$. Without loss of
generality, we only prove that $u_1\in \m^k$ as an example. \\
Choose $y\in \m$ so that $\{y^{\ast}_j, y^{\ast}\}$ form a regular
sequence in $gr_{\m}(R)$ for every $j$. Let $f_j=yf/y_j$; then
$f\in S_k$, therefore there are $u_{ij}\in R$ such that \be
\label{eq4} af_j=\dsum_{i=1}^d u_{ij} z_i.\ee From (\ref{eq3}) and
(\ref{eq4}), we obtain that
$$(u_1y-u_{1j}y_j)z_1\in (z_2, \dots, z_d).$$ Moreover, since
$\{z_1, \dots, z_d\}$ is a regular sequence and $z_i\in \m^{k+1}$
by assumption,
$$(u_1y-u_{1j}y_j)\in \m^{k+1}.$$ It follows that $$u_1\in
(y_j)+\m^k$$ for every $j$ by the choice of $y$. Now one can
conclude that $u_1\in (y_l\cdots y_k)+\m^k$ by the facts that
$\{y^{\ast}_l, y^{\ast}_j\}$ is a regular sequence in $gr_{\m}(R)$
if
$j>l$ for $l=k-1, k-2, \dots$, in particular, $u_1\in \m^k$. \\
\end{proof}
\begin{coll}\label{coll3}
Let $(R, \m)$ be a 2-dimensional regular local ring. Let $k\geq 2$
and ${\bf z}=z_1, z_2$ be regular sequence contained in
$\m^{k+1}$. Then the following hold.
\begin{description}
\item{(i)} $({\bf z}): \m^k\subseteq \m^k$; \item{(ii)} For every
element $a\in ({\bf z}): \m^k$ and every element $f\in \m^k$, if
$af=\sum_{i=1}^2 u_iz_i$ then $u_i\in \m^k$.
\end{description}
\end{coll}
\begin{proof}
(i) follows from the fact that $gr_{\m}(R)$ is Cohen-Macaulay and
$\m^{2k+1}: \m^k= \m^{k+1}\subseteq \m^k$. (ii) follows from the
proof of Theorem~\ref{theo8}.
\end{proof}

Before proving Theorem~\ref{theo7}, we need one more lemma.

\begin{lemma} \label{lem4}
Let $R$ be a Noetherian ring. Let $H$ be an ideal of height $g\geq
2$, $J=(z_1, \dots, z_g)$ an ideal generated by a regular sequence
contained in $H$. Set $I=J: H$. If $I=(a_1, \dots, a_t)$,
$a_iH\subseteq JH$ for every $i$, and $I^2\subseteq J$, then
$I^2=IJ$.
\end{lemma}

\begin{proof}
Since $I^2\subseteq J$, there are $u_{ij}^{(k)}\in R$ such that
$a_ia_j=\sum_{k=1}^g u_{ij}^{(k)} z_k$. To complete the proof it
suffices to show that $u_{ij}^{(k)}\in I$,i.e.,
$u_{ij}^{(k)}H\subseteq J$ for every $k$. For this let $b\in H$;
then by assumption, $$ \ba{rcl}
\sum_{k=1}^g (u_{ij}^{(k)}b)z_k & = & a_ia_jb \\  &  \in  & a_iJH \\
& \subseteq & J^2H \\ & \subseteq & (z_1, \dots, z_g)^2. \ea
$$ Since $z_1, \dots, z_g$ is a regular sequence (although not
necessarily permutable), we see that $u_{ij}^{(k)}b\in (z_1,
\dots, z_g)$ for every $k$.
\end{proof}

\vskip 0.2in ~~\\ {\bf Proof of Theorem~\ref{theo7}:} Let
$I=I_k=J:P^{(k)}$. Then $I\subseteq P^{(k)}$ by
Corollary~\ref{coll1} and  Corollary~\ref{coll3}(i), it follows
that $I^2\subseteq IP^{(k)}\subseteq J$. To complete the proof, by
Lemma~\ref{lem4}, we need only to prove that if $a\in I$, $f\in
P^{(k)}$ and $af=\sum_{i=1}^g u_iz_i$, then $u_i\in P^{(k)}$ for
every $i$. However this follows by Corollary~\ref{coll2} and
Corollary~\ref{coll3}(ii).

\section{The case $t=1$}
We study the links of unmixed radical ideals in a Noetherian ring
in this section. The main result Theorem~\ref{theo9} will
generalize \cite[Theorem~2.3]{cp1} and \cite[Theorem~2.1]{cpv}.
\par Throughout, let $R$ be a Noetherian ring of depth $d\geq g$,
$H$ an unmixed radical ideal of height $g$, and let ${\bf z}=z_1,
\dots, z_g\subseteq H$ be a regular sequence. Set $J=({\bf z})$
and $I=J:H$. Suppose that for every minimal prime $P$ of $H$, one
of the following two conditions holds:\\ $(L_1)$ $R_p$ is not a
regular local ring. \\ $(L_2)$ $R_P$ is a regular local ring of
dimension at least 2 and two elements in the sequence ${\bf
z}=z_1, \dots, z_g$
lie in  the symbolic square $P^{(2)}$. \\
An easy observation is the following.

\begin{lemma} \label{lem5}
$I\subseteq H$ and $I^2\subseteq J$.
\end{lemma}
\begin{proof}
To show $I\subseteq H$, it suffices to show that $I\subseteq P$
for every minimal prime $P$ of $H$. However, if $I$ is not a
subset of $P$ then $P_P=J_P$, which contradicts to $L_1$ or $L_2$.
Thus $I\subseteq H$ is fulfilled. Moreover, $I^2\subseteq
IH\subseteq J$.
\end{proof}

\begin{theo} \label{theo9}
Let $R$ be a Noetherian ring of depth $d\geq g$, $H$ an unmixed
radical ideal of height $g$, and ${\bf z}=z_1, \dots, z_g\subseteq
H$ a regular sequence. Set $J=({\bf z})$ and $I=J:H$. Suppose that
for every minimal prime $P$ of $H$, either $L_1$ or $L_2$ holds.
Then $I^2=JI$.
\end{theo}
\begin{proof}  Write $H=P_1\cap \cdots \cap P_l$. For
$1\leq j\leq l$, let
$$\ba{rcl} T_j & = & \cup_i Ass(R/(z_1, \dots, \hat{z}_i, \dots, z_g))\cap
Spec(R_{P_j})\\ & = &  \cup_i  Min(R/(z_1, \dots, \hat{z}_i,
\dots, z_g))\cap Spec(R_{P_j}) .\ea $$ Since $I$ is not a subset
of every prime in $\cup_j T_j$, we can choose a generating set
$\{a_1, \dots, a_t\}$ of $I$ so that $a_i\notin \cup_j T_j$ for
every $i$. \\ Now, to finish the proof, it is enough to show that
$a_iH\subseteq JH$ for every $i$ by Lemma~\ref{lem4}. For this,
let $b\in H$ and $a$ be any $a_i$; then there are elements $u_j\in
R$ such that \be \label{eq5} ab=\sum_{j=1}^g u_j z_j .\ee To show
$u_j\in H$, we can prove it locally. Without loss of generality,
we will prove that $u_1\in PR_P$ for every minimal prime $P$ of
$H$. \\ Let $b'$ be an arbitrary element in $H$; then there are
elements $u'_j\in R$ such that \be \label{eq6} ab'=\sum_{j=1}^g
u'_j z_j .\ee

~~\\ From (\ref{eq5}) and (\ref{eq6}), we obtain that
$$a(b'u_1-bu'_1)\in (z_2, \dots, z_g).$$ Let $P$ be any minimal prime
of $H$; then $\{a, z_2, \dots, z_g\}$ is a regular sequence in
$R_P$. Since $b'$ is arbitrary, $u_1P_P=u_1H_P\subseteq (b, z_2,
\dots, z_g)R_P$. From the fact that $R_P$ is not regular or at
least one element in $\{z_2, \dots, z_g\}$ is in $P^{(2)}$, we
conclude that $u_1\in PR_P$. \end{proof}

\vskip 0.05in ~~\\  {\bf Acknowledgements}
\par The author would like to thank C. Huneke for helpful
suggestion which improved the exposition of Theorem~\ref{theo6}.

\end{document}